\documentclass[12pt]{amsart}

\usepackage{epsf,amssymb}

\newtheorem{thm}{Theorem}[section]
\newtheorem{prop}[thm]{Proposition}

\newtheorem{lemme}[thm]{Lemma}
\newtheorem{cor}[thm]{Corollary}

\newcommand{\R}{\mathbb{R}}
\newcommand{\Z}{\mathbb{Z}}

\newcommand{\N}{\mathbb{N}}

\newcommand{\bdry}{\partial}
\newcommand{\s}{\vskip.1in}
\newcommand{\n}{\noindent}
\newcommand{\B}{\mathcal{B}}
\newcommand{\F}{\mathcal{F}}

\newcommand{\be}{\begin{enumerate}}
\newcommand{\ee}{\end{enumerate}}

\begin{document}
\title{On the coarse classification of tight contact structures}

\author{Vincent Colin}
\address{Universit\'e de Nantes}
\email{Vincent.Colin@math.univ-nantes.fr}

\author{Emmanuel Giroux}
\address{\'Ecole Normale Sup\'erieure de Lyon}
\email{Emmanuel.GIROUX@umpa.ens-lyon.fr}

\author{Ko Honda}
\address{University of Southern California} 
\email{khonda@math.usc.edu}
\urladdr{http://math.usc.edu/\char126 khonda}

\date{This version: May 30, 2002.}


\begin{abstract}

{We present a sketch of the proof of the following theorems: (1) Every 
$3$-manifold has only finitely many homotopy classes of 2-plane fields 
which carry tight contact structures.  (2) Every closed atoroidal $3$-manifold 
carries finitely many isotopy classes of tight contact structures.} 	

\end{abstract}
\maketitle

In this article we explain how to normalize tight contact structures with 
respect to a fixed triangulation.  Using this technique, we obtain the following 
results:

\begin{thm}  \label{thm: homotopy} 
Let $M$ be a closed, oriented  $3$-manifold.  There are finitely many homotopy 
classes of $2$-plane fields which carry tight contact structures. 
\end{thm}

\begin{thm}  \label{thm: isotopy} 
Every closed, oriented, atoroidal $3$-manifold carries a finite number of tight 
contact structures up to isotopy. 
\end{thm}

P.\ Kronheimer and T.\ Mrowka \cite{KM} had previously shown Theorem~\ref{thm: 
homotopy} for (weakly) symplectically (semi-)fillable contact structures.  Our 
theorem is a genuine improvement of the Kronheimer-Mrowka theorem because there 
exist tight structures which are not fillable~\cite{EH}.

Now, since every Reebless foliation is a limit of tight contact  
structures~\cite{Co4, ET}, we obtain a new proof of a recent result of D.\ Gabai 
\cite{Ga}.

\begin{cor} [Gabai] 
There are finitely many homotopy classes of plane fields which carry Reebless 
foliations. 
\end{cor}

Next, shifting our attention to isotopy classes of contact structures, we see 
that Theorem~\ref{thm: isotopy} complements the following theorem \cite{Co2, 
Co3, HKM}:

\begin{thm}[Colin, Honda-Kazez-Mati\'c]  
Every closed, oriented, irreducible, toroidal $3$-manifold carries infinitely 
many tight contact structures up to isomorphism. 
\end{thm}

Summarizing, we have:
 
\begin{thm} 
A closed, oriented, irreducible $3$-manifold carries infinitely many tight 
contact structures (up to isotopy or up to isomorphism) if and only if it is 
toroidal. 
\end{thm}

\n
A more complete account of the proofs will appear in \cite{CGH}.

\s\n
{\it Acknowledgements.}  This project was started in the fall of 2000 when the 
authors visited the American Institute of Mathematics and Stanford University 
during the Workshop on Contact Geometry.  We wholeheartedly thank the organizers 
Yasha Eliashberg and John Etnyre.   Parts of this work were also realized at 
l'universit\'e de Nantes and the IH\'ES.  Special thanks go to Fran\c cois 
Laudenbach for his careful reading of the paper.

\section{Contact geometry in dimension $3$}

In dimension $3$, there are exactly two types of locally homogeneous 2-plane 
fields $\xi =\ker \alpha$.  The plane field $\xi$ is a {\it contact structure} 
when $\alpha \wedge d\alpha$ is nowhere vanishing, and is a {\it 
foliation} when $\alpha \wedge d\alpha \equiv 0$.  Almost by definition, 
the classification of contact structures and foliations must reflect global 
properties of the ambient manifold $M$.   Contact structures, unlike foliations, 
are very stable objects.  For example, two $C^0$-close 
contact structures are isotopic~\cite{Co1}.  Thus classifying contact structures 
up to isomorphism or isotopy appears to be a reasonable project and has been 
the subject of numerous studies in the past twenty years.

In what follows, we only consider oriented manifolds and oriented contact 
structures which are {\em positive}, i.e., satisfy $\alpha \wedge d\alpha >0$.
Also, denote the metric completion of each component of $A\setminus B$ by 
$\overline{A\setminus B}$.

\subsection{Tight vs.\ overtwisted}

The study of contact $3$-manifolds reveals a dichotomy in the world of 
contact structures:  {\it tight} vs.\ {\it overtwisted} ~\cite{Be, El1, El2}.  A 
contact structure $\xi$ is {\it overtwisted} if there exists an embedded 
$2$-disk $D$ which is tangent to $\xi$ at all its boundary points; otherwise 
$\xi$ is {\it tight}. The boundary of $D$ plays a role analogous to vanishing 
cycles in foliation theory, and there are many analogies between tight contact 
structures and Reebless (or taut) foliations~\cite{ET}.

The classification of overtwisted contact structures up to isotopy coincides 
with the classification of plane fields up to homotopy by the work of  Y.\ 
Eliashberg \cite{El1}; contact structures therefore need tightness to become 
geometrically significant. A similar observation holds in foliation theory: Reeb 
components should be avoided.

\subsection{Convex surfaces}

The main tool for analyzing contact manifolds is the theory of convex surfaces, 
first introduced in~\cite{Gi1}.

Let $S$ be an oriented, properly embedded surface in a contact manifold $(M,\xi 
)$.  We assume that $S$ is either closed or compact with {\it Legendrian} 
boundary (i.e., tangent to $\xi$ at every point).  Then the {\it characteristic 
foliation} $\xi S$ on $S$ is the singular foliation obtained by integrating the 
singular line field $TS \cap \xi$.  It is directed by the orientations of $S$ 
and $\xi$, and is singular precisely when $\xi =TS$.  The singularities are 
generically isolated and of two types: elliptic of index $1$ or hyperbolic of 
index $-1$.  They are positive (or negative) if the  orientation of $\xi$ and 
$TS$ coincide (or not).  The germ of a contact structure near a surface $S$ is 
completely determined by the characteristic foliation $\xi S$.

A surface $S$ is said to be {\it convex} if it is transversal to a vector field 
which preserves $\xi$, i.e., a {\it contact vector field}.  Surprisingly, this 
notion of convexity is generic \cite{Gi1}: every closed surface admits a 
$C^\infty$-small perturbation into a convex surface. Moreover, the convexity can 
be read off from the characteristic foliation $\xi S$; for example, if $\xi S$ 
is singular Morse-Smale, then $S$ is convex \cite{Gi1}.

Let $\gamma$ be a closed embedded Legendrian curve on an embedded surface $S$.  
We define the {\em relative Thurston-Bennequin invariant} $tb(\gamma,S)$ (also 
called the {\em  twisting number} relative to $S$) to be half the number of 
algebraic intersections of $\xi$ and $TS$ along $\gamma$.

If $S$ is a surface with Legendrian boundary and if each component of $\bdry S$ 
has nonpositive relative Thurston-Bennequin invariant, then $S$ can be made 
convex by an isotopy which fixes $\bdry S$, is $C^0$-small in a neighborhood of 
$\bdry S$, and is $C^\infty$-small outside.

Let $S$ be a surface transverse to a contact vector field $X$.  Define the {\it 
dividing set} of $S$ to be $$\Gamma_S =\{ x\in S | X(x)\in 
\xi(x)\}.$$  It is a smooth embedded multicurve transverse to $\xi S$ and 
does not depend on $X$ up to isotopy through multicurves transverse to $\xi S$.  
$\Gamma_S$ inherits a natural orientation from that of $\xi S$.   We say that 
the characteristic foliation $\xi S$ is {\it adapted} to $\Gamma_S$.

The multicurve $\Gamma_S$ captures all of the essential information on $\xi$
in a neighborhood of $S$, according to the following Flexibility Theorem 
\cite{Gi1}:

\begin{thm}[Giroux] 
Let $S$ be a closed surface, $f_0:S\hookrightarrow (M,\xi)$, $g: 
S\hookrightarrow (M',\xi')$ be two embeddings with $f_0(S)$ and $g(S)$ convex, 
and $X$ a contact vector field transverse to $f_0 (S)$.  If the oriented 
multicurves $f_0^{-1} (\Gamma_{f_0 (S)} )$ and $g^{-1} (\Gamma_{g(S)} )$ 
coincide, then there exists an isotopy $(f_t )_{t\in [0,1]}$ transverse to $X$ 
such that $f_1^* (\xi f_1 (S))=g^* (\xi g(S))$.  The same holds if $\bdry S\not 
=\emptyset$ and the characteristic foliations coincide near $\bdry S$; in this 
case we get an isotopy rel $\bdry S$. 
\end{thm}

If $S$ is a surface with Legendrian boundary, an arc $\delta$ of $\Gamma_S$ is 
said {\it boundary-parallel} if one component of $S\setminus \delta$ is a half 
disk which has no other intersections with $\Gamma_S$.

\subsection{Bypasses}

We now introduce the notion of a {\it bypass}, as in \cite{Ho1} and \cite{Et}, 
which allows us to investigate the contact structure outside the neighborhood of 
a convex surface.

Let $D$ be a half-disk embedded in a contact manifold and $\alpha \cup \beta$ be 
the decomposition of $\bdry D$ into two smooth arcs which meet only at their 
endpoints. The disk $D$ is called a {\it bypass} if it satisfies the following: 

\begin{itemize}
\item $\bdry D$ is Legendrian.
\item $D$ is convex, without singularities in its interior.
\item Along $\bdry D$, we have two positive elliptic singularities
at $\bdry \alpha =\bdry \beta$, a positive hyperbolic singularity in the 
interior of $\beta$ and a negative elliptic singularity in the interior of 
$\alpha$. 
\end{itemize}

\n
The Thurston-Bennequin invariant of $\bdry D$ is $-1$, with a contribution of 
$0$ from $\beta$ and $-1$ from $\alpha$.  In other words, $\beta$ is a more 
``efficient'' Legendrian arc than $\alpha$.

\begin{lemme} 
Let $S$ be a convex surface with Legendrian boundary.  Assume there exists a 
boundary-parallel component $\delta$ of $\Gamma_S$ which cuts off a half-disk 
$D_0\subset S$.  (If $S$ is a disk, we assume in addition that $\Gamma_S$ is not 
connected.)  Let $D_1$ be a tubular neighborhood of $D_0$ inside $S$, with $D_1 
\cap \Gamma_S =\delta$.  Denote $\alpha_1 =D_1 \cap \bdry S$.  If $X$ is a 
contact vector field transverse to $S$, then there exists an isotopy of $S$ rel 
$\bdry S$ through surfaces transversal to $X$, which leads to a surface $S'$ 
such that the image $D'\subset S'$ of $D_1$ is a bypass.  
\end{lemme}

The $\alpha_1$ in the lemma plays the role of $\alpha$ in the definition of a 
bypass.  If $S$ has been nicely normalized along $\alpha_1$, the isotopy can be 
chosen to be $C^\infty$-small.

\begin{proof} 
One can draw a foliation on $S$ which is adapted to $\Gamma_S$, and where $D_1$ 
is a bypass. The desired isotopy is then given by the Flexibility Theorem.
\end{proof}

Here is a fundamental result of Y.\ Eliashberg \cite{El2} which classifies tight 
contact structures on the $3$-ball $B^3$.

\begin{thm}\label{theorem: boule} 
Let $\Gamma$ be a connected dividing set on $\bdry B^3$, and $\F$ a foliation 
adapted to $\Gamma$. Then there exists a unique tight contact structure on 
$\B^3$, up to isotopy rel $\bdry B^3$, whose characteristic foliation on $\bdry 
B^3$ is $\F$. 
\end{thm}

\subsection{The generalized Lutz modification}

Let $T\subset (M,\xi)$ be an embedded $2$-torus which meets $\xi$ transversely.
Then there exists a closed tubular neighborhood $U=(\R^2/\Z^2)\times[0,2\pi]$ of 
$T=(\R^2/\Z^2)\times\{\pi\}$ fibered by Legendrian intervals $\{pt\}\times 
[0,2\pi]$.  Then, with respect to the coordinate system $(x,y,t)$ for  $(\R^2 
/\Z^2) \times [0,2\pi]$, $\xi$ is given by: 
$$\cos f(x,y,t)dx -\sin f(x,y,t)dy=0,$$ 
where $f$ is a circle-valued function $(\R^2 /\Z^2) \times [0,2\pi] \rightarrow 
\R/ 2\pi \Z$.  For every $n\in \N$, choose a smooth increasing function $g_n : 
[0,2\pi] \rightarrow \R$ with $g_n (0)=0$ and $g_n (2\pi )=2n\pi$.  The plane 
field defined on $U$ as the kernel of $\alpha_n =\cos (f+g_n )dx-\sin (f+g_n 
)dy$ is a contact structure.  It coincides with $\xi$ on $\bdry U$, and thus can 
be extended by $\xi$ to give a contact structure $\xi_n$ on $M$.  It is called a 
{\it generalized Lutz modification of $\xi$ of index $n$ along $T$}.  (While the 
contact structure $\xi_n$ depends on the choice of tubular neighborhood $U$, 
its isotopy class only depends on $T$.)  

Consider the function $\chi : U\rightarrow [0,1]$ which is $0$ on $\bdry U$ and 
is strictly positive on the interior of $U$. The path of $1$-forms 
$(1-t)\alpha_0 +t\alpha_n +t(1-t)\chi dt$ defines a path of plane fields rel 
$\bdry U$ between $\xi \vert_U$ and $\ker \alpha_n$.
Thus $\xi$ and $\xi_n$ are homotopic through plane fields.

Note that if $(M,\xi)$ is overtwisted, then many generalized Lutz modifications 
do not alter the isotopy class of $\xi$.  For example, if $T$ bounds a solid 
torus and the characteristic foliation $\xi T$ is a suspension, then the new 
contact structure $\xi_n$ is also overtwisted, and is therefore isotopic to 
$\xi$ by Eliashberg's overtwisted classification theorem.  In general, it is not 
clear whether a generalized Lutz modification on an overtwisted contact manifold 
remains overtwisted.  

\s\n
{\bf Exercise:}  Let $(M,\xi)$ be an overtwisted contact manifold and $T\subset 
M$ a pre-Lagrangian incompressible torus.  If $\xi_n$ is obtained from $\xi$ by 
a generalized Lutz modification along $T$, then $(M,\xi_n)$ has a finite cover 
which is overtwisted.

\s\n
{\bf Question:} Is $(M,\xi_n)$ itself overtwisted?

\subsection{Branched surfaces}

\begin{figure} [ht]			
	{\epsfysize=1.5in\centerline{\epsfbox{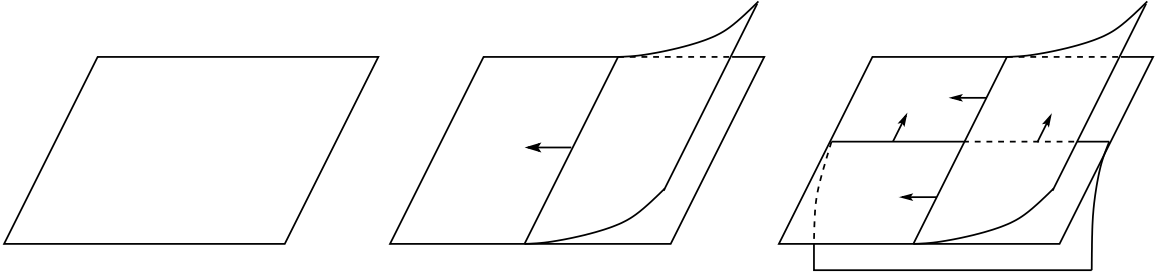}}} 		
	\caption{The three local models for a branched surface.
	The arrows indicate the branching direction.} 				
	\label{branched}
\end{figure}

A {\it branched surface} $\B$ is a topological space such that every $p\in 
\B$ has a neighborhood which is given by one of the three possibilities in 
Figure~\ref{branched}. 
A {\it branched surface with boundary} $\B$ is a topological space locally 
modeled on the three diagrams in Figure~\ref{branched}, plus ``halves'' of the 
diagrams to the left and in the middle, so that $\bdry \B$ is a {\em train 
track}.

\begin{figure} [ht]			
	{\epsfysize=2in\centerline{\epsfbox{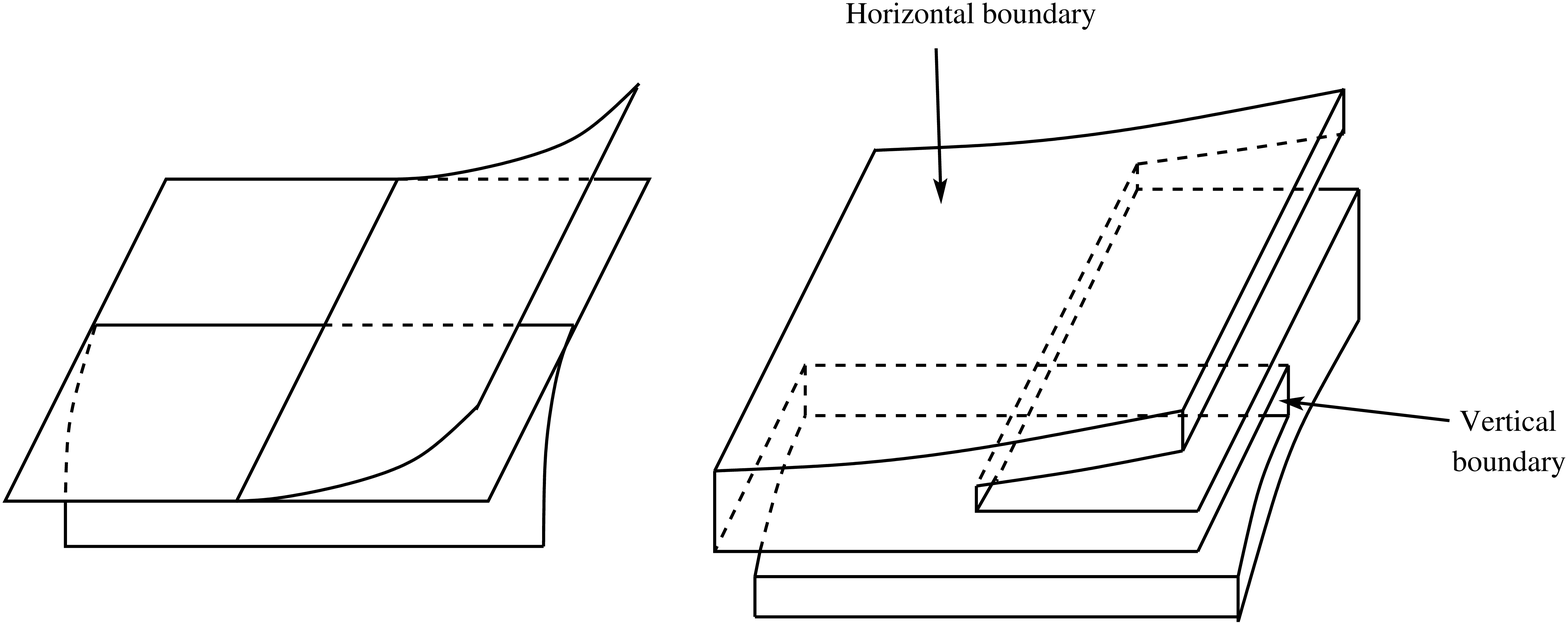}}} 		
	\caption{The fibered neighborhood of a branched surface.} 				
	\label{fibered-nbhd}
\end{figure}

The branch locus $L$ of a branched surface $\B$ is the set of points $p\in \B$ 
for which no neighborhood is modeled on a plane or a half-plane.  It is a 
collection of embedded curves which intersect transversally at double points.  
If $\B$ is embedded inside a $3$-manifold, the neighborhood $N(\B )$ of $\B$ 
admits a fibration by intervals and a projection $\pi : N(\B )\rightarrow \B$ 
(see Figure~\ref{fibered-nbhd}).   The boundary of $N(\B)$ is divided into two 
pieces, the horizontal boundary $\bdry_h N(\B )$ which is transversal to the 
fibration and the vertical boundary $\bdry_v N(\B )$ which is 
tangent to the fibration.

We can now strengthen Theorem~\ref{thm: homotopy}:

\begin{thm}  \label{theorem: Lutz} 
Given a closed $3$-manifold $M$, there exists a finite collection of contact 
structures $\{\xi_1,\dots,\xi_k\}$ and, for each $1\leq i\leq k$, a finite set 
of tori $\mathcal{T}_i=\{T_i^1,\dots,T_i^{k_i}\}$ transverse to $\xi_i$, whose 
union is a branched surface, such that every tight contact structure $\xi$ on 
$M$ is obtained from one of the $\xi_i$ by performing a generalized Lutz 
modification on a subset of $\mathcal{T}_i$. 
\end{thm}

Theorem~\ref{thm: homotopy} follows immediately from Theorem~\ref{theorem: 
Lutz}, since we have already shown that a generalized Lutz modification does not 
alter the homotopy class of the contact structure.

In \cite{Co2, Co3, Gi2, HKM}, generalized Lutz modifications were performed 
along pre-Lagrangian incompressible tori to produce infinitely many 
nonisomorphic tight contact structures. More precisely, if $\xi_n$ is the tight 
contact structure given by the equation $\cos (nt) dx-\sin ( nt) dy=0$ on the 
thickened torus $(\R /\Z )^2 \times [0,2\pi ]$ with coordinates $(x,y,t)$, then 
one can define the {\it torsion} of a contact manifold $(M,\xi )$ as the 
supremum over $n\in \N$ for which there exists a contact embedding $\phi_n :(T^2 
\times [0,2\pi ] ,\xi_n )\hookrightarrow (M,\xi )$.  This invariant of $(M,\xi 
)$, introduced in~\cite{Gi2}, is intended to measure ``how large'' the contact 
manifold is. Notice that if the torsion is finite, then it can be increased by a 
Lutz modification.

In ~\cite{Co2, Co3, HKM}, we essentially prove that every closed, oriented, 
irreducible, and toroidal $3$-manifold carries tight contact structures with 
arbitrarily large and finite torsion.

\section{Normalizing a triangulation in a tight contact manifold}

We now describe the main tool (Proposition~\ref{proposition: branched}) for 
proving Theorems~\ref{theorem: Lutz} and \ref{thm: isotopy}, which echoes the 
theory of normal surfaces of Haken and Kneser~\cite{Ha, Kn}.  This technique is 
also very similar to Gabai's in \cite{Ga}.

Let $\tau$ be a fixed triangulation of a $3$-manifold $M$.  Denote the 
$i$-skeleton of $\tau$ by $\tau^i$.

\subsection{Maximal triangulations}

Given a tight contact structure $\xi$ on $M$, we can isotop $\tau$ so that the 
$1$-skeleton is a Legendrian graph and each face is convex.  Such a 
triangulation will be called a {\em contact triangulation}.  This is easy to 
accomplish because: 

\begin{itemize}

\item every embedded graph can be $C^0$-approximated by a Legendrian graph;
\item the Thurston-Bennequin number $tb(\gamma)$ along each edge $\gamma$,
computed with the trivialization given by any adjacent face, can be made 
strictly negative by a $C^0$-small isotopy of $\tau^1$ rel $\tau^0$;
\item by genericity, after a $C^0$-small isotopy in a neighborhood of 
$\tau^1$ rel $\tau^1$, each face can then be made convex by a $C^\infty$-small 
isotopy rel $\tau^1$. 

\end{itemize}

\n
Note that the triangulation is now singular:  at each vertex $x$, the tangent 
lines to the adjacent edges are all contained in $\xi (x)$.

Denote $TB(\xi ,\tau)=\sum_{F \in \tau^2} tb(\bdry F)$.  If $[\tau]$ is the 
set of all contact triangulations isotopic to $\tau$, then define 
$$TB(\xi,[\tau]) = \max_{\tau'\in[\tau]} TB(\xi,\tau').$$

\begin{figure} [ht]			
	{\epsfysize=2in\centerline{\epsfbox{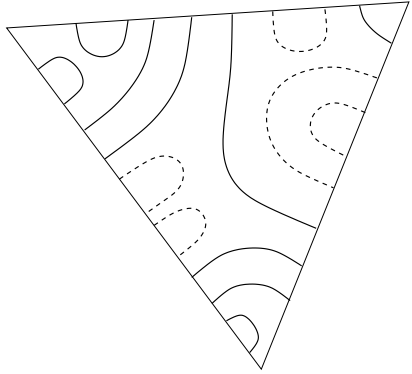}}} 		
	\caption{The dotted dividing curves can be isotoped across, leaving the solid 
	dividing curves.} 					
	\label{face}
\end{figure}

\begin{figure} [ht]			
	{\epsfysize=2in\centerline{\epsfbox{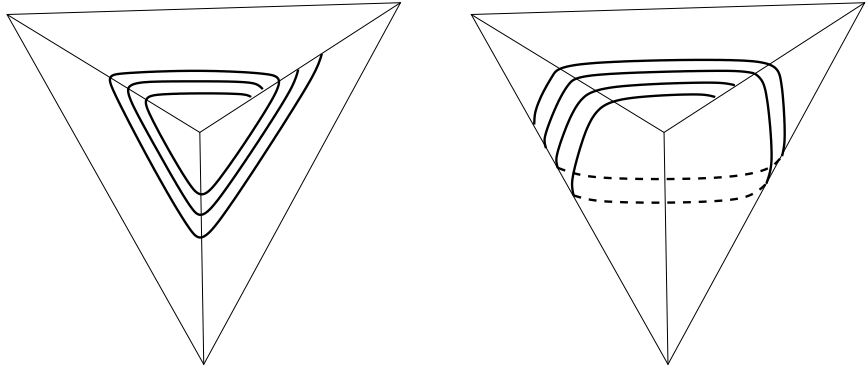}}} 		
	\caption{The ``holonomy'', after the edges are rounded.  In the left-hand 
	diagram, every time the dividing curve wraps around the vertex in a 	
	counterclockwise manner, it moves two units away from the vertex.  The 
	right-hand diagram is similar.} 							
	\label{holonomy} 
\end{figure}

We now assume, after a possible change of notation, that $\tau$ realizes the 
maximum value $TB(\xi,[\tau])$.  Such a triangulation is said to be a {\it 
maximal triangulation} for $\xi$ in the class $[\tau]$.  In that case, we have 
the following:

\begin{lemme} \label{lemma:minimal}   $\mbox{ }$

\be

\item[(a)] For every face $F$ of $\tau$, each component of $\Gamma_F$ is an arc 
whose endpoints are on different edges, except possibly six arcs which are 
boundary-parallel components and whose endpoints are close to a vertex
(see Figure~\ref{face}).

\item [(b)] the ``holonomy'' of the dividing curves around the boundary of a 
$3$-simplex $B$ is as in Figure~\ref{holonomy}. 

\ee
\end{lemme}

\begin{proof} 
We only sketch the proof of (a).  Closed components of $\Gamma_F$ are not 
allowed because the structure is tight.  The existence of a boundary-parallel 
component ``far'' from the vertices would imply the existence of a bypass which 
does not meet the vertices.  Let $\alpha \cup \beta$ be the boundary of the 
bypass, where $\alpha \subset \bdry F$, and isotop the edge $l$ containing 
$\alpha$ to $(l\setminus \alpha) \cup \beta$.  It strictly increases 
$TB(\xi,\tau)$.  See Figure~\ref{simplify}.

\begin{figure} [ht]			
	{\epsfysize=2in\centerline{\epsfbox{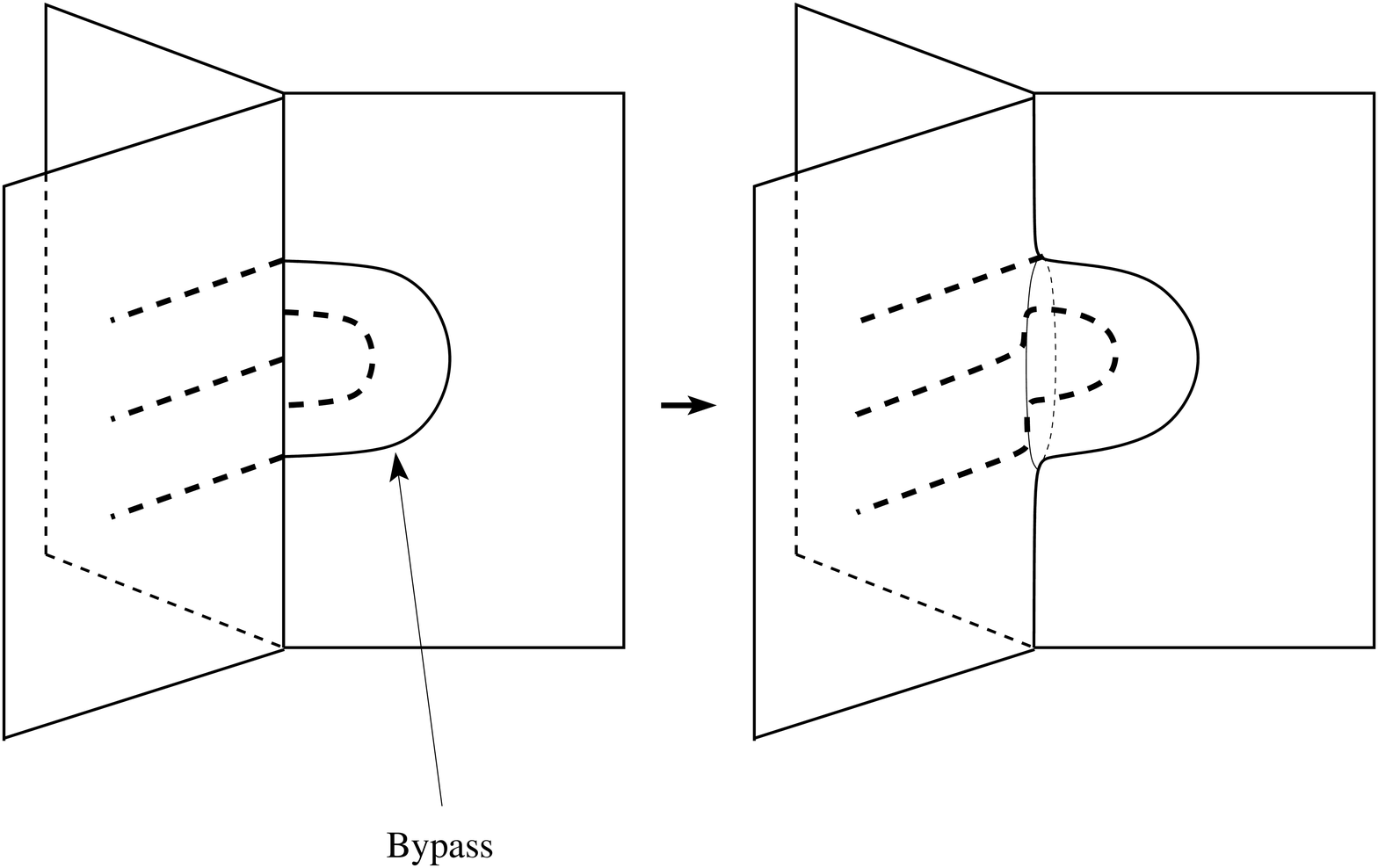}}} 		
	\caption{} 					
	\label{simplify}
\end{figure}

The proof of (b) also relies on the existence of bypasses inside $B$ in case the 
holonomy is not as in Figure~\ref{holonomy} (see~\cite{Ho3}). 
\end{proof}

\subsection{Fibered prisms}

A {\it fibered prism} is a polyhedron $P$ diffeomorphic to $R\times [-1,1]$, 
where $R$ is a triangle or a rectangle, and $\{pt\}\times[-1,1]$ are the fibers.  
The set $R\times \{ -1,1\}$ is called the {\it horizontal boundary}, 
and the set $(\bdry R)\times [-1,1]$ the {\it vertical boundary}.

Let $\xi$ be a tight contact structure on $M$, which has been isotoped so that 
$\tau$ is a maximal triangulation for $\xi$.   Then we partition families of 
parallel arcs of $\Gamma_F$ ($F$ is a face of $\tau$) inside fibered prisms:

\begin{lemme}  \label{lemma: combinatorial}  
For each 3-simplex $B$ of $\tau$, there are at most 5 embeddings $\phi_i :P_i 
\hookrightarrow B$, $i=1,\dots,5$, of fibered prisms $P_i$, such that the 
following are satisfied.

\begin{itemize}

\item The vertical edges (resp.\ faces) of $P_i$ are sent into edges (resp.\ 
faces) of $B$. 

\item $\phi_i (P_i \setminus (\bdry P \times [-1,1]))\subset 
int(B)$. \item $\phi_i (P_i \times \{ -1,1\} )\cap \Gamma_F 
=\emptyset$, for every face $F$ of $B$.

\item The fibered prisms coming from two different 3-simplices intersect along 
rectangles. 

\item For every face $F$ of $B$, at most $C$ components of $\Gamma_F$ are not 
contained in the image of the $\phi_i$'s.   Here $C$ is a universal constant 
which does not depend on $\xi$ or on $\tau$. 

\end{itemize}

\end{lemme}

\begin{figure} [ht]			
	{\epsfysize=2.3in\centerline{\epsfbox{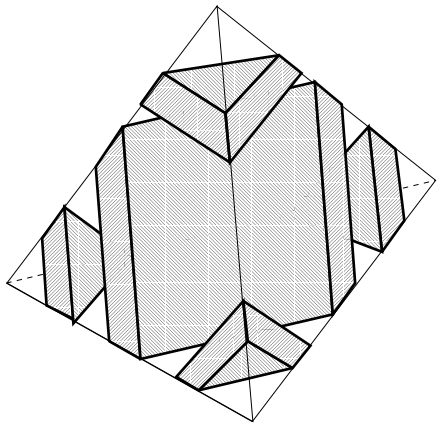}}} 		
	\caption{Fibered prisms.} 							
	\label{prisms} 
\end{figure}

Moreover, by a combination of the Flexibility Theorem and Theorem~\ref{theorem: 
boule}, we can normalize $\xi \vert_{\phi_i (P_i )}$ as in the following lemma, 
after isotoping $\xi$.  (For notational convenience, we will not distinguish 
$P_i$ from its image from now on.)

\begin{lemme}  \label{lemma: polygon}
Fix a nonsingular foliation $\F$ on the horizontal boundary of each prism (the 
same for each face) which coincides with the characteristic foliation near its 
$1$-skeleton.  Then $\xi$ can be isotoped rel $\tau^1$ so that the 
characteristic foliation of $\xi$ equals $\F$ on the horizontal boundary and 
each vertical fiber $\{pt\} \times [-1,1]$ is Legendrian.
\end{lemme}

\subsection{The branched surfaces}

\begin{prop}  \label{proposition: branched}
There exists a finite number of pairs $(\B_1 ,\zeta_1 ),\dots ,(\B_n ,\zeta_n 
)$, where $\B_i$ is a branched surface and $\zeta_i$ a contact structure on 
$M\setminus N(\B_i )$, such that every tight contact structure $\xi$ on $M$, up 
to isotopy, is generated by one of the $(\B_i ,\zeta_i )$, i.e., 
\be
\item $\xi \vert_{M\setminus N(\B_i )} =\zeta_i$;
\item $\xi \vert_{N(\B_i )}$ is tangent to the fibers of $N(\B_i )$.
\ee
Moreover, one may assume that the branched surfaces have empty boundary.
\end{prop}

\begin{proof}  Let $\xi$ be a tight contact structure on $M$, isotoped so that 
$\tau$ is a maximal triangulation for $\xi$.  Apply Lemmas~\ref{lemma: 
combinatorial} and ~\ref{lemma: polygon} to obtain the collection $\{P_1,\dots, 
P_k\}$ of fibered prisms $P_i=R_i\times[-1,1]$ ranging over all the 3-simplices 
$B$ of $\tau$.   For every edge $a\in \tau^1$, denote $a'$ the arc obtained from 
$a$ by deleting a small neighborhood of its endpoints, and $a'\times D^2$ a 
small tubular neighborhood of $a'$ such that the fibration $a'\times \{ pt 
\}$ is Legendrian and coincides with the one given on the $P_i$'s.

\begin{lemme}  \label{lemme: lissage} 
There exists a branched surface $\B$ with boundary, together with a fibered 
neighborhood $N(\B )$ fibered by Legendrian intervals, which is obtained by 
smoothing 
$$K=\left( \bigcup_{1\leq i\leq k} P_i \right) \cup \left(\bigcup_{a\in \tau^1} 
(a'\times D^2)\right).$$  
Moreover, every Legendrian fiber of $N(\B )$ lies inside a Legendrian fiber of 
$K$. \end{lemme}

\begin{proof}[Proof of Lemma~\ref{lemme: lissage}] 
Up to reparametrization, one can always assume that $Q_i \cap Q_j =\emptyset$ 
for $i\neq j$, where we denote $Q_i =R_i \times \{ 0\}$.  Now let $Q_i' 
=\overline{Q_i \setminus \cup_{a\in \tau^1} (a'\times D^2)}$, i.e., we remove 
from $Q_i$ a neighborhood of its vertices.  If $P_i$ intersects $P_j$, we join 
the corresponding edges of $Q_i'$ and $Q_j'$ by a small band transverse
to the fibration.   Denote by $\B'$ the branched surface with boundary 
consisting of the union of the $Q_i'$ and the bands.  Note that $\B'$ does not 
have any triple branch points.

Take a disk $\{pt\} \times D^2 \subset a'\times D^2$ for each edge $a$ of 
$\tau^1$, and squash each point of $\B' \cap (a'\times D^2)$ to the point of 
$\{pt\} \times D^2$ situated on the same fiber.  We obtain a ``singular 
surface'' that can easily be perturbed into a branched surface $\B$; this is 
where we need to introduce triple branch points.  Once the edges of $K$ are 
smoothed, we naturally obtain a fibered neighborhood $N(\B)$ of $\B$ with the 
desired properties. 
\end{proof}

The construction of Lemma~\ref{lemme: lissage} gives rise to a finite number of 
branched surfaces, since each simplex contains at most 5 prisms, and for each 
simplex there are 7 total prism positions to choose from (see 
Figure~\ref{prisms}).  Here we are considering two prisms $P_i$, $P_j$ in a 
3-simplex $B$ to be equivalent if there is an isotopy from $P_i$ to $P_j$ which 
restricts to each edge $a$ to be an isotopy of $P_i\cap a$ to $P_j\cap a$ and on 
each face $F$ to be an isotopy of $P_i\cap F$ to $P_j\cap F$.

We may assume that all the tight contact structures $\xi$ with the same set of 
fibered prisms (and hence the same branched surface $\B$) agree on $\bdry 
N(\B)$, because we can impose the same nonsingular characteristic foliation on 
the horizontal boundaries of the prisms by Lemma~\ref{lemma: polygon}.  
On the complementary regions of $N(\B)$, we have a natural decomposition 
into polyhedra which is inherited from $\tau$.  On the faces of the polyhedra, 
the number of components of the dividing set is universally bounded by 
Lemma~\ref{lemma: combinatorial}.  This leads to a finite number of 
possible dividing curve configurations, and hence to finitely many 
characteristic foliations by the Flexibility Theorem.  But now, according to 
Theorem~\ref{theorem: boule}, a tight contact structure on each polyhedron is 
determined by its restriction to the boundary. 

Now the branched surface $\B$ may have nonempty boundary.  We then use an {\em 
amputation principle} to excise the boundary: if a vertical face of a prism is 
entirely contained in $\bdry_v N(\B)$, then all tight contact structures which 
coincide outside of $N(\B)$ with $\zeta$ and are tangent to the fibers of 
$N(\B)$ coincide on this face.  In particular, they all have the same number of 
components of dividing curves on this face, and the number of dividing curves on 
the other vertical faces of the prism is bounded.  Thus, all the tight contact 
structures which coincide with $\zeta$ on $M\setminus N(\B)$ coincide, up to 
isotopy, with a finite set of models on this prism, by taking into account the 
finite repartition of the parallel dividing curves inside adjacent prisms.  
Therefore, we can amputate this prism from $N(\B)$, at the expense of possibly 
increasing the number of pairs $(\B_i,\zeta_i)$.  This is a finite process which 
must eventually yield a branched surface without boundary (or an empty branched 
surface).  \end{proof}

\section{Tight contact structures carried by the fibered neighborhood of a 
branched surface}

We now analyze contact structures which are {\it carried} by the fibered 
neighborhood of a branched surface $N(\B )$, i.e., are tangent to the fibers and 
are fixed on the boundary.  

Let $\mathcal{S}$ be the set of tight contact structures tangent to the fibers 
of $N(\B )$ (with a given germ along $\bdry N(\B)$), and let $B_1,\dots,B_d$ be 
the components of $\overline{\B\setminus L}$ (recall $L$ is the singular locus 
of $\B$).  If we fix $\xi_0\in \mathcal{S}$, then for any given $\xi\in 
\mathcal{S}$ one can define a weight system $w_\xi$ on $\B$: $w_\xi(B_i )$, 
$i=1,\dots,d$, is the difference of rotation between $\xi$ and $\xi_0$ along a 
fiber $\pi^{-1}(p)$, where $p\in B_i$.  It is a set of integers which does not 
depend on the choice of $p\in B_i$.  Note that, unlike the relative 
Thurston-Bennequin invariant, a left twist contributes positively to the 
weight.  It is easy to verify that $\xi \in \mathcal{S}$ is determined up to 
isotopy by its weight $w_\xi = (w_\xi(B_1), \dots, w_\xi(B_d))\in \Z^d$.  
Also, each $w_\xi(B_i)$ is bounded below by the positive contact structure 
condition.  By amputating the $B_i$ for which $w_\xi(B_i)$ is negative, we may 
assume that $w_\xi(B_i)$ is nonnegative for each $\xi\in \mathcal{S}$. 

To each smooth edge $A$ of $L$, we associate a linear equation 
$$\phi_A(x_1,\dots,x_d)=x_i -(x_j+x_k)=0,$$ 
with $(x_1,\dots,x_d)\in\Z^d$, where $B_i$, $B_j$ and $B_k$ are the components 
of $\B\setminus L$ which meet along $A$, and $B_j$ and $B_k$ branch out of 
$B_i$.  Let $V$ be the subspace of $\R^d$ of solutions to the system 
$\{\phi_A=0 \mbox{ $|$ $\forall$ smooth edges $A$ of $L$}\}$.  Each weight $w$ 
is an element of $V\cap \Z_{\geq 0}^d$, where $\Z_{\geq 
0}^d=\{(x_1,\dots,x_d)\in \Z^d| x_i\geq 0 \mbox{ } \forall i\}$.

Define a partial order $\leq$ on $\Z^d$ by:
$$(x_1 ,\dots ,x_d)\leq (x_1' ,\dots x_d' )\; {\rm if}\; x_i \leq x_i' \;{\rm for}\;
i=1,\dots ,d.$$

\n
There is only a finite number of minimal elements of $V\cap \Z_{\geq 0}^d$ with 
respect to $\leq$, and we denote them $u_1,\dots ,u_k$.  They generate the 
monoid $V\cap \Z_{\geq 0}^d$.  On the other hand, we have the following standard 
fact:

\begin{lemme} 
The isotopy classes of embedded surfaces carried by $N(\B )$ (i.e., transverse 
to the fibers) are in bijection with elements of $V\cap \Z_{\geq 0}^d$.
\end{lemme}

\n
{\bf Remark.}  Every surface carried by $N(\B)$ is transverse to a contact 
structure and is thus either a torus or a Klein bottle.

\s\n
Let $\{T_1 ,\dots,T_k\}$ be the collection of surfaces carried by $N(\B)$ and 
determined by the weights $u_1 ,\dots ,u_k$. The following lemma, together with 
Proposition~\ref{proposition: branched} and the remark below, yields 
Theorem~\ref{theorem: Lutz}.

\s\n
Assume all the $T_i$'s are tori.   Then we have:

\begin{lemme} 
Every contact structure $\xi \in \mathcal{S}$ is obtained from $\xi_0$ by a 
collection of generalized Lutz modifications along $T_1,\dots,T_k$. 
\end{lemme}

\begin{proof} 
Given $\xi \in \mathcal{S}$, we can write $w_\xi=\sum_{1\leq i\leq k} n_i u_i$, 
$n_i \in \Z_{\geq 0}$.  The contact structure $\xi'$ obtained from $\xi_0$ by 
generalized Lutz modifications of index $n_i$ along $T_i$ is tangent to the 
fibers of $N(\B)$ and has the same weight as $\xi$.  Therefore $\xi$ and 
$\xi'$ are isotopic rel $\bdry N(\B)$. 
\end{proof}

\n
{\bf Remark.}  If $T_i$ is a  Klein bottle, then we can replace it by the torus 
$T_i'=\bdry N(T_{i})$, where $N(T_i)$ is the tubular neighborhood of $T_i$.  
Denote $T_i'=T_i$ if $T_i$ is a torus.  Then the set $\mathcal{S}$ is still 
generated by Lutz modifications along the $T_i'$, but starting from a finite set 
of structures (not only $\xi_0$), obtained from $\xi_0$ by performing 
$\pi$-twists along the Klein bottles.

\section{Isotopy and beyond}

The proof of Theorem~\ref{thm: isotopy} follows from a more in-depth analysis of 
tight contact structures carried by the fibered neighborhood of these branched 
surfaces.  We use the classification of tight contact structures on the 
thickened torus (see~\cite{Gi4, Ho1}) to build a process which allows us to 
decrease the weights in case the manifold is atoroidal.  In the general case, 
one should be able to prove that there exists a finite number of tight contact 
structures of bounded torsion up to isomorphism.  Moreover, these techniques 
apply to manifolds $M$ with boundary, provided $\bdry M$ is not a torus and we 
prescribe a fixed characteristic foliation along $\bdry M$.  We then obtain  
relative versions of Theorems~\ref{thm: homotopy} and  \ref{thm: isotopy}.
Finally, we have the following:

\begin{thm} Given an oriented topological knot type $\mathcal{K}$ in the 
standard tight contact 3-sphere $(S^3,\xi_{std})$ and an integer $n$, there 
exists, up to contact isotopy, a finite number of Legendrian knots in 
$\mathcal{K}$ with Thurston-Bennequin invariant equal to $n$. \end{thm}

\end{document}